\documentclass[11pt, reqno]{amsart}
\usepackage{amssymb}
\usepackage{amscd}
\usepackage{verbatim,ifthen}
\usepackage{color}
\usepackage{latexsym}
\usepackage{tikz}
\usepackage{tikz-cd}
\usepackage{mathrsfs}
\usepackage{wrapfig}
\usetikzlibrary{shapes}
\usepackage{color}
\usetikzlibrary{arrows.meta}
\usepackage{bbm}
\usetikzlibrary{matrix}
\usetikzlibrary{calc}
\usetikzlibrary{arrows,intersections}
\usepackage{pgfplots}
\usepackage{multicol}
\usepackage{array}
\newcolumntype{M}[1]{>{\centering\arraybackslash}m{#1}}

\usepackage[colorlinks, linkcolor=black, citecolor=blue, linktocpage,urlcolor=blue,backref=page]{hyperref}
\renewcommand*{\backref}[1]{}
\renewcommand*{\backrefalt}[4]{[{\tiny%
    \ifcase #1 Not cited.%
          \or Cited on page~#2.%
          \else Cited on pages #2.%
    \fi%
    }]}
\addtocontents{toc}{\setcounter{tocdepth}{1}}

\usepackage{amsmath}

\numberwithin{equation}{section}

\let\oldtocsection=\tocsection
 
\let\oldtocsubsection=\tocsubsection

\renewcommand{\tocsection}[2]{\hspace{0em}\oldtocsection{#1}{#2}}
\renewcommand{\tocsubsection}[2]{\hspace{1em}\oldtocsubsection{#1}{#2}}

\def\XXint#1#2#3{{\setbox0=\hbox{$#1{#2#3}{\int}$ }
\vcenter{\hbox{$#2#3$ }}\kern-.6\wd0}}

\usepackage{eucal}
\usepackage{calc}  
\usepackage{enumitem} 
\usepackage{tensor}
\usepackage{graphicx,wrapfig,lipsum}
\usepackage{etoolbox}
\usepackage{marginnote}
\usepackage{lipsum}
\makeatletter
\patchcmd{\@mn@margintest}{\@tempswafalse}{\@tempswatrue}{}{}
\patchcmd{\@mn@margintest}{\@tempswafalse}{\@tempswatrue}{}{}
\reversemarginpar 
\makeatother
\usepackage{scrextend}

\makeatletter
\DeclareRobustCommand\widecheck[1]{{\mathpalette\@widecheck{#1}}}
\def\@widecheck#1#2{%
    \setbox\z@\hbox{\m@th$#1#2$}%
    \setbox\tw@\hbox{\m@th$#1%
       \widehat{%
          \vrule\@width\z@\@height\ht\z@
          \vrule\@height\z@\@width\wd\z@}$}%
    \dp\tw@-\ht\z@
    \@tempdima\ht\z@ \advance\@tempdima2\ht\tw@ \divide\@tempdima\thr@@
    \setbox\tw@\hbox{%
       \raise\@tempdima\hbox{\scalebox{1}[-1]{\lower\@tempdima\box
\tw@}}}%
    {\ooalign{\box\tw@ \cr \box\z@}}}
\makeatother
\title{Some remarks on the Wu--Yau Theorem}
\author{Kyle Broder}
\address{The University of Queensland,  St. Lucia,  QLD 4067, Australia}
\email{k.broder@uq.edu.au}
\thanks{The author was supported by the Australian Government through the Australian Research Council's Discovery Projects funding scheme (project DP220102530). }
\usepackage[utf8]{inputenc}
\usepackage[T1]{fontenc}
\begin{document}

\maketitle

\begin{abstract}
The Wu--Yau theorem asserts that a compact K\"ahler manifold with negative holomorphic sectional curvature admits a cohomologous metric with negative Ricci curvature.  We introduce a conjectural positive analog of the Wu--Yau theorem and present several examples verifying the conjecture.  We also produce an abundance of examples of K\"ahler--Einstein metrics on K\"ahler surfaces of general type that do not have negative holomorphic sectional curvature.  This begins to address a question raised by Diverio.
\end{abstract}

\section{Introduction}
A compact complex manifold $X$ is  \textit{hyperbolic} if every holomorphic map $\mathbf{C} \to X$ is constant (see,  e.g., \cite{KobayashiHyperbolicComplexSpaces}).   One of the main folklore conjectures (generalizing a conjecture of Kobayashi \cite{KobayashiHyperbolicComplexSpaces}) is that a compact hyperbolic manifold should have ample canonical bundle. In particular,  they are all projective, canonically polarized, and thus admit a unique invariant K\"ahler--Einstein metric \cite{Aubin,Yau1976}.   The Ahlfors--Schwarz lemma \cite{Ahlfors} implies that a compact K\"ahler manifold $(X, \omega)$ with negative holomorphic sectional curvature $\text{HSC}(\omega)<0$ is hyperbolic.  While this class of manifolds does not exhaust the class of compact hyperbolic manifolds \cite{Demailly}, such manifolds form an important and general class of hyperbolic manifolds.  For such manifolds,  we have the following (see \cite{WuYau1,TosattiYang} and the references therein):

\subsection*{Theorem 1.1}\label{WuYauTheorem}
(Wu--Yau theorem).  Let $(X, \omega)$ be a compact K\"ahler manifold with $\text{HSC}(\omega)<0$.  Then $X$ is projective and canonically polarized.  In particular,  there is a smooth function $\varphi \in \mathcal{C}^{\infty}(X,\mathbf{R})$ such that \begin{eqnarray}\label{WUYAU1}
\text{HSC}(\omega) \ < \ 0 \ \implies \ \text{Ric}(\omega_{\varphi}) \ < \ 0, \hspace*{1cm} \omega_{\varphi} : = \omega + dd^c \varphi.
\end{eqnarray} 

\hfill

Given the importance and attention that the above result has received (see e.g., \cite{Diverio,BroderSBC,BroderSBC2,BroderStanfield}),  it is natural to ask (c.f., \cite{Diverio}) whether \eqref{WUYAU1} holds with the signs reversed:

\subsection*{Question 1.2}\label{Question1.2}
Let $(X, \omega)$ be a compact K\"ahler manifold.  If $\text{HSC}(\omega)>0$, does there exist a (cohomologous) metric $\omega_{\varphi} : = \omega + dd^c \varphi$ such that $\text{Ric}(\omega_{\varphi})>0$? \\

Hirzebruch surfaces $\mathcal{F}_n : = \mathbf{P}(\mathcal{O}_{\mathbf{P}^1} \oplus \mathcal{O}_{\mathbf{P}^1}(n))$ admit K\"ahler metrics $\omega$ with $\text{HSC}(\omega)>0$ for all $n \in \mathbf{N}_0$, but  for $n>1$, $\mathcal{F}_n$ does not admit any K\"ahler metric with positive Ricci curvature \cite{Hitchin}.  In particular,  not only can one not produce a cohomologous metric of positive Ricci curvature from a K\"ahler metric of positive holomorphic sectional curvature,  there may be no metrics of positive Ricci curvature at all! The `positive analog' of \nameref{WuYauTheorem} is quite subtle, and a conjectural picture has remained out of reach.  By considering the details of \cite{Hitchin}  and building a heuristic picture of the behavior of $\text{HSC}(\omega)$ and $\text{Ric}(\omega)$ under bimeromorphic modifications (e.g., blow-ups), we propose the following conjecture:

\subsection*{Conjecture 1.3}\label{ConjecturePositive}
Let $(X, \omega)$ be a compact K\"ahler manifold with $\text{Ric}(\omega)>0$.  Then there is a cohomologous metric $\omega_{\varphi} : = \omega + dd^c \varphi$ such that $\text{HSC}(\omega_{\varphi})>0$. \\

Yang \cite{YangHSCYau} has shown that a compact K\"ahler manifold $(X, \omega)$ with $\text{HSC}(\omega)>0$ is projective and rationally connected (i.e., any two points lie in the image of a rational curve $\mathbf{P}^1 \to X$).  \nameref{ConjecturePositive} can therefore be interpreted as an analytic formulation of the well known result of Campana \cite{CampanaFano} and Koll\'ar--Miyaoka--Mori \cite{KollarMiyaokaMori} that Fano manifolds (i.e., K\"ahler manifolds with ample anti-canonical bundle) are rationally connected.  To the author's knowledge,  no analytic proof of this result is known. 

Let us remark that the requirement that the metric of positive holomorphic sectional curvature in \nameref{ConjecturePositive} is cohomologous to the given metric of positive Ricci curvature is motivated only by \nameref{WuYauTheorem}.  It would be interesting to understand whether there are compact K\"ahler manifolds with positive Ricci curvature that admit K\"ahler metrics with positive holomorphic sectional curvature but are not cohomologous.

The following result produces a rich number of examples that verify \nameref{ConjecturePositive}:

\subsection*{Theorem 1.4}
\nameref{ConjecturePositive} holds for the following manifolds: \begin{itemize}
\item[(1)] K\"ahler C--spaces whose irreducible factors are either: \begin{itemize}
\item[(i)] The classical K\"ahler C--spaces of type $A,B,C, D$. 
\item[(ii)] The exceptional K\"ahler C--spaces of type $(E_6,\alpha_p)$ for $p= 1, ..., 6$,  $(E_7, \alpha_1)$,  $(E_7, \alpha_2)$, $(E_7, \alpha_6)$, $(E_7, \alpha_7)$, $(E_8, \alpha_1)$, $(E_8, \alpha_8)$, $(F_4, \alpha_1)$, $(F_4, \alpha_4)$, and $(G_2, \alpha_2)$. 
\end{itemize}
\item[(2)] Fano manifolds $X$ that are the total space of a surjective holomorphic submersion $p : X \to S$ with fibers $X_s : = p^{-1}(s)$ and base $S$ admitting K\"ahler metrics of positive holomorphic sectional curvature.
\end{itemize}

Statement (1) and (2) are essentially due to Itoh \cite{Itoh} who showed that the invariant K\"ahler--Einstein metrics on the irreducible K\"ahler C--spaces have positive holomorphic sectional curvature.  In particular, all Fano Hermitian symmetric spaces verify \nameref{ConjecturePositive}.  Statement (2) follows from the recent work of Chaturvedi--Heier \cite{ChaturvediHeier}.  For the reader's convenience,  we will give a brief reminder of the Fano K\"ahler C--spaces in $\S 2$.  The details of the curvature calculation are left to \cite{Itoh}, where they are already well-presented.

Formulating \nameref{ConjecturePositive} as a `positive analog' of \nameref{WuYauTheorem} may be fruitful in addressing some long-standing problems concerning the holomorphic sectional curvature.  For instance,  a resolution of \nameref{ConjecturePositive} for Fano surfaces would address the following long-standing question raised by Yau: Does the blow-up of $\mathbf{P}^2$ at two points admit a K\"ahler metric with $\text{HSC}>0$?  

Gabriel Khan pointed out to me that there are some circumstances where positive holomorphic sectional curvature may imply positive Ricci curvature.  For instance, if the holomorphic sectional curvature is positive and suitably pinched, the Ricci curvature will also be positive.  In fact, if the holomorphic sectional curvature is positive and suitably pinched,  the orthogonal bisectional curvature is positive, and by \cite{Chen}, the underlying manifold is biholomorphic to $\mathbf{P}^n$.

Diverio \cite{Diverio} raised the question of whether the cohomological perturbation in \eqref{WUYAU1} is non-trivial, in general.  While we certainly expect the perturbation to be non-trivial in general, it would be of interest if one could classify or produce obstructions to a compact complex manifold admitting a metric with both negative holomorphic sectional curvature and negative Ricci curvature.  Further, note that as a consequence of the Aubin--Yau theorem \cite{Aubin,Yau1976},  we have the secondary implication: \begin{eqnarray}\label{SecondImplication}
\text{HSC}(\omega) \ < \ 0 \ \implies \ \text{Ric}(\omega_{\varphi}) \ < \ 0 \ \implies \ \text{Ric}(\omega_{\text{KE}}) \ = \ - \omega_\text{KE},
\end{eqnarray} where $\omega_{\text{KE}} = \omega_{\varphi} + dd^c u$ is the unique invariant K\"ahler--Einstein metric on $X$.

Motivated by Diverio's question, it is natural to ask:

\subsection*{Question 1.5}
Can one classify the (compact) K\"ahler--Einstein manifolds whose Einstein metric has negative holomorphic sectional curvature? \\

It was brought to my attention by Gabriel Khan that the same question was posed in \cite{KhanZheng}.  He also informed me that Tong \cite{TongPinching} showed that the K\"ahler--Ricci flow, starting from a K\"ahler metric whose holomorphic sectional curvature is negatively pinched, converges to a K\"ahler--Einstein metric that is uniformly equivalent to the initial metric. In particular, if the holomorphic sectional curvature of the initial metric is sufficiently pinched,  then it remains so for all time, and thus yields a K\"ahler--Einstein metric with negative holomorphic sectional curvature. Obstructions to the K\"ahler--Einstein metric on bounded pseudoconvex domains in $\mathbf{C}^n$ having negatively pinched holomorphic sectional curvature were given by Cho \cite{Cho}.

Stefano Trapani brought to my attention the very recent paper \cite{Sarem}, where an example of a compact K\"ahler manifold with a K\"ahler metric of negative Ricci curvature and quasi-negative (i.e., non-positive everywhere and negative at a point) holomorphic sectional curvature was produced, but the manifold is not even Kobayashi hyperbolic! In particular, the manifold has no (Hermitian) metric of negative (Chern) holomorphic sectional curvature.

From the developments concerning the `geography problem' for surfaces of general type, we have the following:

\subsection*{Theorem 1.6}
Let $(X, \omega)$ be a K\"ahler--Einstein surface with $\text{Ric}(\omega) = - \omega$ and $\text{HSC}(\omega)<0$. Then \begin{eqnarray}\label{CheungChernClass}
c_2(X) \ \leq \ 3c_1^2(X).
\end{eqnarray} In particular,  the following examples do not have K\"ahler--Einstein metrics with negative holomorphic sectional curvature: Barlow surfaces, Burniat surfaces,  Campadelli surfaces, Catanese surfaces,  Godeaux surfaces,  Horikawa surfaces,  Keum--Naie surfaces,  Oliverio surfaces, Todorov surfaces. \\

The Chern class inequality \eqref{CheungChernClass} is due to Cheung \cite{CheungArticle}, but no examples were provided.  Hence, the main contribution of the above theorem is to exhibit a vast number of examples that fail to have K\"ahler--Einstein metrics with negative holomorphic sectional curvature.  Cheung also provides a sufficient condition for the K\"ahler--Einstein metric to have negative $\text{HSC}<0$.  This is expressed in terms of the Chern--Weil representatives of the Chern classes, and not the Chern classes themselves.  In particular,  verifying Cheung's criterion is significantly more difficult, since these representatives are seldom given explicitly.   Cheung makes some effort to produce a class of examples that satisfy his criterion.  We will show in \nameref{BochnerCurvatureTensor} however, that this class of manifolds is of no interest.

To conclude the introduction,  let us mention that the asymmetry between the conditions $\text{HSC}_{\omega}>0$ and $\text{HSC}_{\omega}<0$ can be seen very explicitly by considering the Gauduchon connections (see, e.g., \cite{BroderThesis, BroderStanfield}). In \cite[Corollary 4.14]{BroderStanfield}, a monotonicity theorem for the Gauduchon holomorphic sectional curvature was discovered.  This shows that, at least for Hermitian non-K\"ahler metrics, ${}^c \text{HSC}_{\omega}<0$ is the strongest condition, while ${}^c \text{HSC}_{\omega}>0$ is the weakest condition, among the  Gauduchon holomorphic sectional curvatures. This offers a candidate explanation for why ${}^c \text{HSC}_{\omega} \leq - \kappa_0 < 0$ implies Kobayashi hyperbolicity  without any completeness or K\"ahler assumption on the Hermitian metric. On the other hand, ${}^c \text{HSC}_{\omega}>0$ does not imply anything close to rationally connectedness -- the standard metric on the Hopf surface $\mathbf{S}^1 \times \mathbf{S}^3$ has ${}^c \text{HSC}_{\omega}>0$, but there are no rational curves at all!

\subsection*{Acknowledgements}
Many of the ideas and examples presented in this work have originated from lectures I delivered between 2020 and 2023. I was eventually given sufficient encouragement to write such things down.  I would like to thank Nigel Hitchin, Stefano Trapani,  Simone Diverio,  Christoph B\"ohm,  Jorge Lauret,   Peter Petersen,  Ramiro Lafuente,  Finnur L\'arusson, Franc Forstneri$\check{\text{c}}$,  Filippo Bracci, Wolfgang Ziller,  and James Stanfield for their interest in this work and valuable discussions.  Part of this work began during my Ph.D.  under the supervision of Ben Andrews and Gang Tian.  I am sincerely grateful for their unwavering support and encouragement.

After posting the first version of this manuscript to the arXiv, I received some very valuable comments from Gunnar \TH\'or Magn\'usson, Gabriel Khan, and Stefano Trapani.  I am thankful for these comments that have greatly improved the manuscript.

\section{A Positive Analog of the Wu--Yau Theorem}
Recall that a \textit{rational curve} (respectively, \textit{entire curve}) is a non-constant holomorphic map $\mathbf{P}^1 \to X$ (respectively, $\mathbf{C} \to X$).  The Ahlfors Schwarz lemma \cite{Ahlfors} implies that a compact K\"ahler manifold $(X, \omega)$ with $\text{HSC}(\omega) < 0$ has no rational or entire curves.   The holomorphic sectional curvature $\text{HSC}(\omega)$ of a K\"ahler metric is similar in strength to the Ricci curvature $\text{Ric}(\omega)$ in the sense that they both dominate the scalar curvature \cite{BergerHBC} (for the Ricci curvature, this is obvious), and both are dominated by the holomorphic bisectional curvature $\text{HBC}(\omega)$ (for the holomorphic sectional curvature, this is obvious).   The relationship between $\text{HSC}(\omega)$ and $\text{Ric}(\omega)$ is particularly curious, however:

\subsection*{Example 2.1}
By adjunction, a smooth hypersurface $X_d \subset \mathbf{P}^n$ of degree $d > n+1$ has ample canonical bundle. In particular, by the Aubin--Yau theorem \cite{Aubin, Yau1976}, $X_d$ admits a K\"ahler(--Einstein) metric of negative Ricci curvature. There are Fermat hypersurfaces of degree $d > n+1$ that contain complex lines, however,  and therefore, cannot have $\text{HSC}(\omega)<0$ (see, e.g., \cite{WuYau1}).  Another example of a surface of general type that is not hyperbolic was constructed by Hirzebruch \cite{Hirzebruch} (see also \cite{CheungArticle}). 

On the other hand,  the Wu--Yau theorem \cite{ WuYau1,TosattiYang} tells us that a compact K\"ahler manifold with $\text{HSC}(\omega)<0$ is projective with ample canonical bundle.  Hence, again by the Aubin--Yau theorem \cite{Aubin, Yau1976},  for a compact K\"ahler manifold $(X, \omega)$ with $\text{HSC}(\omega)<0$, there is a smooth function $\varphi \in \mathcal{C}^{\infty}(X, \mathbf{R})$ such that $\text{Ric}(\omega_{\varphi}) < 0$, where $\omega_{\varphi} : = \omega + dd^c \varphi$. 

Given the significant interest and activity concerning the Wu--Yau theorem, it has been of interest to understand whether a `positive analog' of the Wu--Yau theorem holds.  Naively,  it is natural to ask whether compact K\"ahler manifolds $(X, \omega)$ with $\text{HSC}(\omega)>0$ admit a cohomologous metric $\omega_{\varphi} : = \omega + dd^c \varphi$ with $\text{Ric}(\omega_{\varphi})>0$ (i.e., \nameref{Question1.2}).

\subsection*{Example 2.2}
Hitchin \cite{Hitchin} showed that all Hirzebruch surfaces $$\mathcal{F}_n \ : = \  \mathbf{P}(\mathcal{O}_{\mathbf{P}^1} \oplus \mathcal{O}_{\mathbf{P}^1}(n)), \hspace{1cm} n \in \mathbf{N}_0,$$ admit Hodge metrics of $\text{HSC}>0$.  These metrics are given by a warped product construction.  For $n>1$, however,  the first Chern class $c_1(\mathcal{F}_n)$ is not positive, and hence, there is no K\"ahler metric with $\text{Ric}>0$.

%\subsection*{Explain the Details of Hirzebruch's Proof}
\nameref{ConjecturePositive} emerges from two heuristics concerning the behavior of the holomorphic sectional curvature and Ricci curvature under bimeromorphic modifications (e.g., blow-ups). The first observation is the following basic fact: Let $\Phi : \widetilde{X} \to X$ be a bimeromorphic modification of $X$ (e.g., the blow-up of $X$ at a point).  The exceptional divisor of $\Phi$ is isomorphic to $\mathbf{P}^{n-1}$, and hence, $\widetilde{X}$ supports a number of rational curves $\mathbf{P}^1 \to \widetilde{X}$.  In particular, $\widetilde{X}$ never admits a metric of negative holomorphic sectional curvature. In this sense, we have the following heuristic: 

\subsection*{Heuristic 1}
Bimeromorphic modifications (in particular, blow-ups) tend to `increase' the holomorphic sectional curvature.  \\

It should be emphasized that it is far from known whether the (total space $\widetilde{X}$ of a) blow-up $\mu : \widetilde{X} \to X$ of a compact K\"ahler manifold $X$ with positive holomorphic sectional curvature has positive holomorphic sectional curvature.  In fact,  Magn\'usson \cite{Magnusson} has shown that the curvature of any metric of the form $\mu^{\ast} \omega_X + t \alpha$, for $\alpha$ a Hermitian form that looks like the Fubini--Study metric near the exceptional divisor, always has negative holomorphic sectional curvature in some directions, for sufficiently small $t>0$. In particular, the natural approach to producing a metric of positive holomorphic sectional curvature on $\widetilde{X}$ is doomed to fail.

To understand the corresponding behavior of the Ricci curvature,  it suffices to understand the behavior of the first Chern class $c_1 : = c_1(K_X^{-1})$.  Here, we restrict our attention to the case of complex surfaces, since we get a precise picture, following \cite[p. 68]{Hitchin}: The Nakai--Moishezon criterion tells us that  $c_1 >0$ if and only if $c_1^2 >0$ and $c_1 \cdot [ C ] >0$, for any curve $C$. Let $\Phi : \widetilde{Y} \to Y$ be the blow-up of $Y$ of a smooth K\"ahler surface $Y$. For a curve $C$ in $Y$,  we have \begin{eqnarray*}
c_1(Y) \cdot [C] \ = \ \Phi^{\ast} c_1 \cdot \Phi^{\ast}[C]  \ = \ c_1(\widetilde{Y}) \cdot \Phi^{\ast} [C] \ = \ c_1(\widetilde{Y}) \cdot ([\widetilde{C}] + m [E]),
\end{eqnarray*}

where $E$ is the exceptional divisor of $\Phi$, and $\widetilde{C} : = \Phi^{-1}(C)$.  Further, $c_1^2(Y) = c_1^2(\widetilde{Y}) +1$.  In particular,  $c_1(\widetilde{Y}) >0$ implies $c_1(Y)>0$.  In the same sense as Heuristic 1,  we see that:

\subsection*{Heuristic 2}
Bimeromorphic modifications (in particular, blow-ups) tend to `decrease' the Ricci curvature. \\

The only relations between the holomorphic sectional curvature and the Ricci curvature that maintain Heuristic 1 and 2 are the Wu--Yau theorem, \nameref{WuYauTheorem} and the conjectural `positive analog' we proposed before, namely \nameref{ConjecturePositive}.

\subsection*{A Brief Reminder of K\"ahler C--spaces}
The natural testing ground for \nameref{ConjecturePositive} is the class of Fano manifolds with symmetries.  By a theorem of Kobayashi \cite{Kobayashi1961Annals},  a compact K\"ahler manifold with positive Ricci curvature (i.e., a Fano manifold) is simply connected.  Hence,  by a theorem of Matsushima \cite{Matsushima3}, a homogeneous Fano manifold admits an invariant K\"ahler--Einstein metric that is unique up to scaling.   A simply connected compact K\"ahler manifold that is homogeneous (in the sense that the holomorphic isometry group acts transitively) is called a \textit{K\"ahler C--space}.  A K\"ahler C--space $X$ is said to be \textit{irreducible} if it cannot be written as $X \simeq Y \times Z$ for positive-dimensional complex manifolds $Y,Z$.  Observe that if $X$ is a Fano K\"ahler C--space, then $Y,Z$ are Fano K\"ahler C--spaces. Moreover, by the result of Chaturvedi--Heier \cite{ChaturvediHeier} if $Y, Z$ have positive holomorphic sectional curvature, then $X$ has positive holomorphic sectional curvature. Hence, it suffices to consider irreducible Fano K\"ahler C--spaces.

For an irreducible K\"ahler C--space,  the group of holomorphic isometries is a simply connected, simple complex Lie group. We write $G$ for this group and denote by $\mathfrak{g}$ its Lie algebra and by $\mathfrak{h}$ a Cartan subalgebra.  Set $\ell : = \dim_{\mathbf{C}} \mathfrak{h}$ and write $\Delta$ for the root system of $\mathfrak{g}$ with respect to $\mathfrak{h}$. By fixing a fundamental root system $\{ \alpha_1, ..., \alpha_{\ell} \}$ of $\Delta$, we have an ordering of the root system $\Delta = \Delta^+ \cup \Delta^-$, and a decomposition of the Lie algebra $$\mathfrak{g} \ = \ \mathfrak{h} \oplus \bigoplus_{\alpha \in \Delta} \mathfrak{g}_{\alpha}.$$ 

For $1 \leq r \leq \ell$, set $\Delta_r^+(k) : = \left \{ \sum_i n_i \alpha_i \in \Delta^+ : n_r = k \right \}$, and $\Delta_r^+ : = \bigcup_{k >0} \Delta_r^+(k)$. Let $P$ be the subgroup whose Lie algebra is $$\mathfrak{p} \ = \ \mathfrak{h} \oplus \bigoplus_{\alpha \in \Delta - \Delta_r^+} \mathfrak{g}_{\alpha},$$ then $P$ is a parabolic subgroup (i.e., $M = G / P$ is compact). This yields a K\"ahler C--space with $b_2=1$. Conversely, any K\"ahler C--space with $b_2=1$ is given in this way.  We denote this space by $(\mathfrak{g}, \alpha_r)$. Itoh \cite{Itoh} computed the curvature of the invariant K\"ahler--Einstein metric on K\"ahler C--spaces. He showed that the invariant metric has positive holomorphic sectional curvature on all K\"ahler C--spaces for which $\Delta_r^+(k) = \emptyset$ for $k \geq 3$. This condition is satisfied by all four classical sequences $A, B, C, D$, for all $r$, and for the exceptions $(E_6, \alpha_p)$ for $1 \leq p \leq 6$,  $(E_7, \alpha_p)$ for $p \in \{ 1,2,6,7 \}$, $(E_8, \alpha_p)$ for $p \in \{ 1, 8 \}$,  $(F_4, \alpha_1)$, $(F_4, \alpha_4)$, and $(G_2, \alpha_2)$.

\section{Surfaces Whose K\"ahler--Einstein Metric Does Not Have Negative Holomorphic Sectional Curvature}
We briefly remind the reader of the proof of Cheung's theorem \cite{CheungArticle} that a compact K\"ahler surface with a K\"ahler--Einstein metric of negative holomorphic sectional curvature satisfies $c_2(X) \leq 3c_1^2(X)$.  Let $(X, \omega)$ be a compact K\"ahler surface with $\text{Ric}(\omega) = \lambda \omega$, for $\lambda <0$. Fix a point $p \in X$ and let $\{ e_1, e_2 \}$ be a local frame for $T_p^{1,0} X$.  The holomorphic sectional curvature $\text{HSC}_{\omega}$ in the unit direction $v = (v_1, v_2) \in T_p^{1,0}X$ is given by $$\text{HSC}_{\omega}(v) \ = \ \sum_{i,j,k, \ell=1}^2 R_{i \overline{j} k \overline{\ell}} v_i \overline{v}_j v_k \overline{v}_{\ell}.$$ Computing the partial derivatives of $\text{HSC}_{\omega}(v)$ with respect to the real and imaginary parts of $v_1, v_2$, if the minimum occurs of $\text{HSC}_{\omega}$ occurs in the direction $e_1$,  then the only non-zero components of the curvature $R$ are those with precisely two indices indices equal to $1$ and precisely two indices equal to $2$.  Write $H_{\min} : = R_{1 \overline{1} 1 \overline{1}} = \text{HSC}_{\omega}(e_1)$.

The holomorphic sectional curvature in this distinguished frame is then \cite[(3.1)]{CheungArticle}: \begin{eqnarray*}
\text{HSC}_{\omega}(v) &=& H_{\min} + 2(2R_{1 \overline{1} 2 \overline{2}} - H_{\min}) | v_1 \overline{v}_2 |^2 + 2 \text{Re} \left( R_{1 \overline{2} 1 \overline{2}} (v_1 \overline{v}_2)^2 \right).
\end{eqnarray*}

From \cite[(3.4)]{CheungArticle},  the holomorphic sectional curvature is negative at $p \in X$ if and only if \begin{eqnarray}\label{(3.4)}
H_{\min} + \frac{1}{2}(2 R_{1 \overline{1} 2 \overline{2}} - H_{\min} + | R_{1 \overline{2} 1 \overline{2}} | ) & < & 0.
\end{eqnarray} 

Let $\gamma_1$, $\gamma_2$ denote the Chern--Weil functions of the curvature tensor of the K\"ahler--Einstein metric $\omega$. That is, the Chern classes are given by $$c_1^2(X) \ = \ \int_X \gamma_1^2 \omega^2, \hspace*{1cm} c_2(X) \ = \ \int_X \gamma_2 \omega^2.$$ In terms of the curvature, we have \begin{eqnarray*}
\gamma_1 \ = \  H_{\min} + R_{1 \overline{1} 2\overline{2}},  \hspace*{1cm}  \gamma_2 \ = \ \frac{1}{2} (H_{\min}^2 + 2 (R_{1 \overline{1} 2 \overline{2}})^2 + | R_{1 \overline{2} 1 \overline{2}} |^2).
\end{eqnarray*}

From \eqref{(3.4)},  the K\"ahler--Einstein metric has negative holomorphic sectional curvature at $p \in X$ if $\gamma_2(p) < \gamma_1^2(p)$ (see \cite[Lemma 3.9]{CheungArticle}). Moreover,  by \cite[Theorem 3.12]{CheungArticle}: If a compact K\"ahler surface has an Einstein metric of negative holomorphic sectional curvature, then $c_2(X)  \leq  3 c_1^2(X)$. 

\subsection*{Example 3.1}
Cheung provides no examples of surfaces satisfying or not satisfying this condition, except for ball quotients (i.e.,  when $c_1^2(X) = 3 c_2(X)$).  Using the number of developments concerning the geography problem of complex surfaces, we exhibit several examples: \begin{itemize}
\item[(1)] Barlow surfaces \cite{Barlow}: A simply connected surface of general type with $p_g=0$ that is homeomorphic but not diffeomorphic to $\mathbf{P}^2$ blown up at 8 points. These surfaces have $c_1^2 = 1$ and $c_2 = 11$.
\item[(2)] Burniat surfaces \cite{Burniat}: Surfaces $S$ of general type with $p_g =0$ and $K_S^2 = 2, 3, 4, 5, 6$ constructed as $\mathbf{Z}_2 \oplus \mathbf{Z}_2$ Galois covers of $\mathbf{P}^2$ branched over certain configurations of nine lines (see also \cite[p.  301--302]{BHPV}).  These surfaces have $c_1^2 = 2$ and $c_2 = 10$.
\item[(3)] Campadelli surfaces \cite{Campadelli,ReidBook}: A minimal surface $S$ of general type with $p_g=q=0$ and $K_S^2=2$.  These surfaces have $c_1^2 = 2$ and $c_2=10$.
\item[(4)] Catanese surfaces \cite{Catanese}: A minimal surface $S$ of general type with $K_S^2=2$.  These surfaces have the same Hodge numbers as the Campadelli surfaces and have $c_1^2 = 2$ and $c_2 = 10$.
\item[(5)] Godeaux surfaces \cite{Godeaux,ReidBook}: A minimal surface $S$ of general type with $p_g = q =0$ and $K_S^2 =1$.  These surfaces have $c_1^2  = 1$ and $c_2 = 11$.
\item[(6)] Horikawa surfaces \cite{Horikawa1, Horikawa2}: Minimal surfaces $S$ of general type with $K_S^2 = 2(p_g-2)$ or $K_S^2 = 2p_g - 3$.  These surfaces lie on, or immediately above the Noether line.  It can be shown that all Horikawa surfaces have $c_1^2/c_2 < 1/3$.
\item[(7)] Keum--Naie surfaces \cite{Keum,Naie}: Families of minimal surfaces of general type with $K_X^2 = 4$ and $p_g=0$ constructed as double covers of an Enriques surface with eight nodes.  These surfaces have $c_1^2 = 1$ and $c_2 = 11$.
\item[(8)] Oliverio surfaces \cite{Oliverio}: Surfaces $S$ of general type with $K_S^2 = 8$, $p_g =4$, and $q=0$ whose canonical system is basepoint free.  Their canonical models are the general complete intersections of bidegree $(6,6)$ in the weighted projective space $\mathbf{P}(1^2, 2, 3^2)$.  Oliverio surfaces have $c_1^2 = 8$ and $c_2 = 52$.
\item[(9)] Todorov surfaces \cite{Todorov}: Minimal surfaces $S$ of general type with $p_g = 1$, $q=0$ and $2 \leq K_S^2 \leq 8$.  Some Todorov surfaces have $c_1^2/c_2 \geq 1/3$,  but many have $c_1^2/ c_2 < 1/3$. 
\end{itemize}

\subsection*{Remark 3.2}\label{BochnerCurvatureTensor}
Cheung also provides a sufficient condition for the K\"ahler--Einstein metric to have negative holomorphic sectional curvature (see, \cite[Lemma 3.9]{CheungArticle}).  In contrast to the necessary condition,  the sufficient condition is expressed in terms of the Chern--Weil functions, not the Chern classes.  This makes it particularly difficult to check.  From \cite[p. 110]{CheungArticle}, if the Chern--Weil function $\gamma_2$ is constant,  then the inequality $\gamma_2 < \gamma_1^2$ is expressed in terms of the Chern classes.  In \cite[Corollary 3.10]{CheungArticle}, Cheung asserts that this is true if the Bochner curvature tensor has constant length. Again,  no examples of such manifolds are given.  Observe,  however, that if the Bochner curvature tensor has constant length,  it is, in particular, parallel. Matsumoto--Shukichi \cite{MatsumotoShukichi} showed that a K\"ahler manifold with parallel Bochner curvature tensor is either locally symmetric or has vanishing Bochner curvature tensor (such manifolds are said to be Bochner--K\"ahler).  By an old result of Kamishima \cite{Kamishima}, the only compact Bochner--K\"ahler manifolds are the compact quotients of the known symmetric ones.

\end{document}